\documentclass[12pt]{amsart}
\usepackage{amscd,amssymb}
\usepackage[graph,frame,poly,arc]{xy}  
\usepackage[plainpages,backref,urlcolor=blue]{hyperref}

\theoremstyle{plain}
\newtheorem{thm}[subsection]{Theorem}

\newtheorem{prop}[subsection]{Proposition}
\newtheorem{cor}[subsection]{Corollary}

\theoremstyle{definition}
\newtheorem{rk}[subsection]{Remark}
\newtheorem{definition}[subsection]{Definition}
\newtheorem{ex}[subsection]{Example}

\newtheorem{question}[subsection]{Question}

\numberwithin{equation}{section}
\newcommand{\OO}{{\mathcal O}}

\newcommand{\A}{{\mathcal A}}
\newcommand{\B}{{\mathcal B}}

\newcommand{\CC}{{\mathcal C}}

\newcommand{\al}{{\alpha}}
\newcommand{\be}{{\beta}}

\newcommand{\Z}{\mathbb{Z}}

\newcommand{\C}{\mathbb{C}}
\newcommand{\PP}{\mathbb{P}}

\DeclareMathOperator{\genus}{genus}

\DeclareMathOperator{\indeg}{indeg}
\DeclareMathOperator{\Bour}{Bour}

\begin{document}

\title [A new hierarchy for complex plane curves]
{A new hierarchy for complex plane curves}

\author[Takuro Abe]{Takuro Abe}
\address{Department of Mathematics, 
Rikkyo University,
Tokyo 171-8501, Japan.}
\email{abetaku@rikkyo.ac.jp}

\author[Alexandru Dimca]{Alexandru Dimca}
\address{Universit\'e C\^ ote d'Azur, CNRS, LJAD and INRIA, France and Simion Stoilow Institute of Mathematics,
P.O. Box 1-764, RO-014700 Bucharest, Romania.}
\email{Alexandru.Dimca@univ-cotedazur.fr}

\author[Piotr Pokora]{Piotr Pokora}
\address{Department of Mathematics,
University of the National Education Commission Krakow,
Podchor\c a\.zych 2,
PL-30-084 Krak\'ow, Poland.}
\email{piotrpkr@gmail.com, piotr.pokora@uken.krakow.pl}



\begin{abstract} 
We define the type of a plane curve as the initial degree of the corresponding Bourbaki ideal. Then we show that this invariant behaves well with respect to the union of curves.
Curves of type $0$ are precisely the free curves, while curves of type $1$ are the plus-one generated curves. In this paper, we first show that line arrangements and conic-line arrangements can exhibit all the theoretically possible types. In the second part, we study the properties of the curves of type $2$ and construct families of line arrangements and conic-line arrangements of this type.
\end{abstract}
 
\maketitle

\section{Introduction} 

Let $S=\C[x,y,z]$ be the polynomial ring in three variables $x,y,z$ with complex coefficients, and let $C:f=0$ be a reduced curve of degree $d\geq 3$ in the complex projective plane $\PP^2$. The ring $S$ and all the graded $S$-modules considered below are graded with respect to the natural grading.
We denote by $J_f$ the Jacobian ideal of $f$, i.e. the homogeneous ideal in $S$ spanned by the partial derivatives $f_x,f_y,f_z$ of $f$, and by $M(f)=S/J_f$ we denote the corresponding graded quotient ring, called the Jacobian (or Milnor) algebra of $f$.
Consider the graded $S$-module of Jacobian syzygies of $f$, namely
$$D_0(f)=\{(a,b,c) \in S^3 \ : \ af_x+bf_y+cf_z=0\}.$$
According to Hilbert's Syzygy Theorem, the graded Jacobian algebra $M(f)$ has a minimal resolution of the form
\begin{equation}
\label{res1}
0 \to F_3 \to F_2 \to F_1 \to F_0,
\end{equation}
where clearly $F_0=S$, $F_1=S^3(1-d)$ and the morphism $F_1 \to F_0$ is given by
$$(a,b,c) \mapsto af_x + bf_y + cf_z.$$
With this notation, the graded $S$-module of Jacobian syzygies $D_0(f)$ has the following minimal resolution
$$0 \to F_3(d-1) \to F_2(d-1).$$
We say that $C:f=0$ is an {\it $m$-syzygy curve} if the module  $F_2$ has rank $m$. Then the module $D_0(f)$ is generated by $m$ homogeneous syzygies, say $r_1,r_2, \ldots ,r_m$, of degrees $d_j=\deg r_j$ ordered such that $$1 \leq d_1\leq d_2 \leq \ldots \leq d_m.$$ 
We call these degrees as the {\it exponents} of the curve $C$, and $r_1, \ldots ,r_m$ as a {\it minimal set of generators } for the module $D_0(f)$. Note that $d_1=mdr(f)$ is the minimal degree of a non-trivial Jacobian relation in $D_0(f)$. It is known that
\begin{equation}
\label{boundG}
d_1\leq d_2 \leq d_3 \leq d-1,
\end{equation}
see for instance \cite[Theorem 2.4]{3syz}.
On the other hand,  for line arrangements we have a better bound, namely 
\begin{equation}
\label{bound}
d_{m} \leq d-2 ,
\end{equation}
see \cite[Corollary 3.5]{Sch}, and for conic-line arrangement we have
\begin{equation}
\label{boundQ}
d_{m} \leq d-1,
\end{equation}
see \cite[Corollary 1.5]{CMreg}.

Let $I_f$ denote the saturation of the Jacobian ideal $J_f$ with respect to the maximal ideal ${\bf m}=(x,y,z)$ in $S$. Consider the following local cohomology group, usually called the Jacobian module of $f$, 
 $$N(f)=I_f/J_f=H^0_{\bf m}(M(f)).$$
We set $n(f)_k=\dim N(f)_k$ for any integer $k$ and we introduce the {\it freeness defect of the curve} $C$ by the formula
$$\nu(C)=\max _j \{n(f)_j\}$$ as in \cite{3syz}.
Note that $C$ is free if and only if $N(f)=0$ and hence $\nu(C)=0$, and $C$ is nearly free if and only if $\nu(C)=1$, see for instance \cite{Dmax, DStRIMS} for the definitions of free and nearly free curves.
Recall that the total Tjurina number of a given curve $C \subset \mathbb{P}^{2}$  is defined to be
$$\tau(C) = \sum_{p \in {\rm Sing}(C)} \tau_{p},$$ 
where ${\rm Sing}(C)$ denotes the set of all singular points of $C$, and $\tau_p$ is the Tjurina number of the singularity $(C,p)$.

Consider the general form of the minimal resolution for the Milnor algebra $M(f)$ of a curve $C:f=0$ that is assumed to be not free, namely
\begin{equation}
\label{res2A}
0 \to \oplus_{i=1} ^{m-2}S(-c_i) \to \oplus_{j=1} ^mS(1-d-d_j)\to S^3(1-d)  \to S,
\end{equation}
with $c_1\leq \ldots \leq c_{m-2}$ and $d_1\leq \ldots \leq d_m$.
It follows from \cite[Lemma 1.1]{HS} that one has
\begin{equation}
\label{res2B}
c_j=d+d_{j+2}-1+\varepsilon_j,
\end{equation}
for $j=1,\ldots,m-2$ and some integers $\varepsilon_j\geq 1$. Using \cite[Formula (13)]{HS}, it follows that one has
\begin{equation}
\label{res2C}
d_1+d_2=d-1+\sum_{j=1} ^{m-2}\varepsilon_j.
\end{equation}
This formula implies the following result, see \cite[Theorem 2.3]{3syz} for details.

\begin{thm}
\label{thmPO1}
Let $C:f=0$ be a reduced plane curve of degree $d$ and let $d_1$ and $d_2$ be the minimal degrees of a minimal system of generators for the module of Jacobian syzygies $D_0(f)$ as above.
Then the following conditions hold:
\begin{enumerate}

\item The curve $C$ is free if and only if $d_1+d_2=d-1$.

\item The curve $C$ is plus-one generated if and only if $d_1+d_2=d$.

\item The curve $C$ is neither free nor plus-one generated if and only if $d_{1}+d_{2} > d$.

\end{enumerate}

\end{thm}
For the definition of plus-one generated curves, see \cite{Abe18,3syz}.
In view of this result, it seems natural to introduce the following invariant, which leads to an interesting hierarchy among all plane curves.

\begin{definition}
\label{deftypek}
Let $C:f=0$ be a reduced plane curve of degree $d$ and let $d_1$ and $d_2$ be the two smallest degrees of a minimal system of generators for the module of Jacobian syzygies $D_0(f)$ as above.
We then say that $C$ has type $t(C)$, where $t(C)$ is a positive integer defined by the formula
$$t(C)=d_1+d_2+1-d.$$
\end{definition}

The fact that this invariant plays a key role in the study of plane curves
follows from the following alternative definition. Let $C:f=0$ be a reduced curve in $\PP^2$ and $r_1$ be a minimal degree generator of the graded $S$-module $D_0(f)$. Then the corresponding Bourbaki ideal $B(C,r_1) \subset S$ 
of $D_0(f)$ was introduced in 
\cite[Section 5]{Split}, and \cite [Section 3]{3syz}, after being considered
 in a non-explicit way by du Plessis and Wall in \cite{duPCTC}, and the first author in \cite{Dmax}. Its properties were further explored in the recent paper \cite{JNS}. It follows from the definition of the Bourbaki ideal $B(C,r_1)$ that the type of the curve $C$ is just the initial degree of the
homogeneous ideal $B(C,r_1)$, namely
\begin{equation}
\label{eqT}
t(C)=\indeg B(C,r_1)= \min \{s \in \Z \ : \ B(C,r_1)_s \ne 0 \}.
\end{equation}
Recall that the Bourbaki degree $\Bour(C)$ of the curve $C$ is defined by
\begin{equation}
\label{eqBou1}
\Bour(C)=\deg B(C,r_1)
\end{equation}
and it satisfies the equality
\begin{equation}
\label{eqBou2}
\Bour(C)=(d-1)^2-d_1(d-d_1-1)-\tau(C),
\end{equation}
see \cite[Theorem 5.1 (1)]{Split} and \cite[Theorem 2.1]{JNS}.
It follows that, if $d_1 \leq (d-1)/2$, then
\begin{equation}
\label{eqBou3}
\Bour(C)=\nu(C),
\end{equation}
see Theorem \ref{deff} below.

Note that a curve $C$ is free exactly when $t(C)=0$ and also if and only if $\Bour(C)=0$.
On the other hand, a curve $C$ is plus-one generated exactly when $t(C)=1$. A plus-one generated curve is a 3-syzygy curve with exponents $(d_1,d_2,d_3)$ satisfying $d_1+d_2=d$, see \cite[Theorem 2.3 (2)]{3syz}. For such a curve one has
\begin{equation}
\label{eqTau10}
\tau(C)=(d-1)^2-d_1(d-d_1-1)+d_2-d_3-1,
\end{equation}
see \cite[Proposition 2.1 (4)]{3syz}. It follows that, for a plus-one generated curve $C$ with exponents $(d_1,d_2,d_3)$ one has
\begin{equation}
\label{eqTau11}
\Bour(C)=d_3-d_2+1,
\end{equation}
and hence $\Bour(C)=1$ if and only if $d_2=d_3$. This occurs when $C$ is a nearly free curve, which is a special type of plus-one generated curve. 

The type of a curve is a positive integer bounded by $2d-2$ in view of \eqref{boundG}, and this upper bound is attained, for instance, by a smooth curve. For line arrangements the corresponding upper bound is
$2d-4$ in view of \eqref{bound}, see Corollaries \ref{cor10} and \ref{cor11} below for more on this subject.

Our first result says that the type behaves well with respect to the union of two curves.  Consider two reduced curves $C_1:f_1=0$ and $C_2:f_2=0$, with no common irreducible component. Then one has the following easy result, concerning the union curve $C=C_1 \cup C_2: f=f_1 \cdot f_2=0$, see \cite[Theorem 5.1]{DIS}.

\begin{prop}
\label{propA}
With the above notation, one has
$$\max\{mdr(f_1),mdr(f_2)\} \leq mdr(f) \leq \min \{mdr(f_1)+\deg C_2, mdr(f_2)+ \deg(C_1)\}.$$
\end{prop}
The type of curves behaves similarly with respect to the union of curves,
namely we have the following result.

\begin{thm}
\label{thm1}
Let $C_1:f_1=0$ and  $C_2:f_2=0$ be two curves in $\PP^2$ with no irreducible component in common. Then the union curve $C=C_1 \cup C_2:f=f_1 \cdot f_2=0$ satisfies
$$t(C) \leq \min\{t(C_1)+\deg(C_2), t(C_2)+\deg(C_1)\}.$$
\end{thm}
\begin{rk}
If $C_1$ is free and $C_2$ is a line, then this result is already known -- see \cite{POG,MP}.
\end{rk}
\begin{cor}
\label{cor1}
If $C$ is a free curve and $Q$ is a conic intersecting $C$ in a finite number of points, then the union $C'=C \cup Q$ is a curve of type at most $2$.
Similarly, if $C$ is a plus-one generated curve and $L$ is a line, which is not an irreducible component of $C$, then  the union $C'=C \cup L$ is a curve of type at most $2$.
\end{cor}

The next result tells us what happens for the type $t(C)$, in most cases, when the curve $C_2$ is generic,
namely smooth and transversal to $C_1$.

\begin{thm}
\label{thm10}
Let $C_1:f_1=0$ be a curve and let $d_1\leq d_2$ be the two smallest degrees of a minimal system of generators for the module of Jacobian syzygies $D_0(f_1)$.
Let $C_2:f_2=0$ be a  smooth curve such that the intersection  $C_1 \cap C_2$ consists of $\deg C_1 \cdot \deg C_2$ points. If
$$d_2 \leq  \deg C_1-2,$$
then the union curve $C=C_1 \cup C_2:f=f_1 \cdot f_2=0$ satisfies
 $$t(C)=t(C_1)+ \deg C_2.$$
 More precisely, if $d_1'\leq d_2'$ denote the first two smallest degrees of a minimal system of generators for the module of Jacobian syzygies $D_0(f)$, then 
 $$d_1'=d_1+ \deg C_2  \text{ and } d_2' =d_2+ \deg(C_2) \leq \deg(C)-2.$$

\end{thm}

\begin{cor}
\label{cor10}
Let $\A$ be a line arrangement  in $\PP^2$ containing at least 3 lines. Then one has
$$0 \leq t(\A) \leq \deg(\A)-3.$$
Conversely, for any integers $t $ and $d $ such that $0 \leq  t \leq d-3$,
there are line arrangements  in $\PP^2$ such that $t(\A)=t$ and $\deg(\A)=d$.
\end{cor}

Proving a result similar to Corollary \ref{cor10} for conic-line arrangements is more complicated. The smooth conic $C_1$ in Example \ref{ex10} satisfies $t(C_1)=d-1$, while the conic-line arrangement $C$ from the same example satisfies $t(C)=d-2$. It is easy to check that 2 smooth conics meeting in 4 points give rise to conic-line arrangement $\CC$ with $t(\CC)=2=d-2$, since $d=\deg(\CC)=4$.
We believe that such examples are exceptional, and ask the following question. 

\begin{question}
\label{q10}
Is it possible to list all the conic-line arrangements $\CC$ satisfying
$$ t(\CC) \geq deg(\CC)-2  \  ?$$
\end{question}
On the positive side, we have the following analog of Corollary \ref{cor10} for conic-line satisfying an additional condition.

\begin{cor}
\label{cor11}
Let $\CC:f=0$ be a conic-line arrangement  and let $d_1\leq d_2$ be the two smallest  degrees of a minimal system of generators for the module of Jacobian syzygies $D_0(f_1)$.
If $d_2 \leq deg( \CC) -2$, then
$$0 \leq t(\CC) \leq \deg(\CC)-3.$$
Conversely, for any integers $t$ and $d$ such that $0 \leq t \leq d-3$,
there are conic-line arrangements $\CC$ containing at least one conic in $\PP^2$ such that $t(\CC)=t$ and $\deg(\CC)=d$.

In particular, for any degree $d \geq 3$, there are free conic-line arrangements $\CC$ containing at least one conic  such that  $\deg(\CC)=d$.
\end{cor}

Before presenting further results, we would like to make the following observation regarding the type.
\begin{rk}
\label{rk10A}
The type in the class of line arrangements is not determined by the combinatorics.
Indeed, the pair of arrangements 
$$\A: f=xy(x-y-z)(x-y+z)(2x+y-2z)(x+3y-3z)(3x+2y+3z)$$
$$(x+5y+5z)(7x-4y-z)=0$$
and
$$\A': f'=xy(4x-5y-5z)(x-y+z)(16x+13y-20z)(x+3y-3z)(3x+2y+3z)$$$$(x+5y+5z)(7x-4y-z)=0$$
constructed by Ziegler \cite{Z} and Yuzvinski \cite{Y} have the same combinatorics, but the exponents are
$(5,6,6,6)$ and $(6,6,6,6,6,6)$, respectively. Therefore, the type of $\A$ is $3$, while the type of $\A'$ is $4$. 
\end{rk}

One can make a loose analogy of the type of curves introduced in Definition \ref{deftypek} with Arnold's notion of modality of function germs. In the case of function germs, the $0$-modal (alias the simple singularities), the $1$-modal (alias the unimodal singularities) and the $2$-modal singularities are classified, see for instance \cite{AGV}, and the complexity of the corresponding geometry, e.g. reflected in their Milnor lattices, is increasing as the modality increases.

Following this path, from now on in this note on we focus on curves $C$ having type $2$. The exponents and the corresponding Tjurina numbers of the curves of type $2$ are described in the following result, where the notation from \eqref{res2C} is used.

\begin{prop}
\label{propB}
Let $C$ be a curve of type $2$.
Then one of the following two cases occurs.

\begin{enumerate}

\item $C$ is a 3-syzygy curve and $\varepsilon_1=2$. If $(d_1,d_2,d_3)$ are the exponents of $C$, then
$$\tau(C)=d_1^2+d_1d_2+d_2^2-2d_1-2d_2-2d_3.$$

\item $C$ is a 4-syzygy curve and $\varepsilon_1= \varepsilon_2= 1$. If $(d_1,d_2,d_3,d_4)$ are the exponents of $C$, then
$$\tau(C)=d_1^2+d_1d_2+d_2^2-2d_1-2d_2-d_3-d_4+1.$$

\end{enumerate}
Conversely, if $m=3$ and $\varepsilon_1=2$ or $m=4$ and $\varepsilon_1=\varepsilon_2=1$, then $t(C)=2$.

\end{prop}
In case $(1)$ (resp. $(2)$), {\it we say that $C$ is a curve of type} $2A$ (resp. $2B$).
For the curves $C$ of type 2, we also have simple formulas for the freeness defect $\nu(C)$, as shown below.

\begin{prop}
\label{propCA}
Suppose that $C$ has type $2A$ and exponents $(d_1,d_2,d_3)$.
Then the following hold.
\begin{enumerate}

\item If $d_1 <d_2-2$, then
$\nu(C)=2(d_3-d_2)+4 \geq 4.$
\item If $d_1=d_2-2$, then
$\nu(C) = 2(d_{3}-d_{1}) \geq 4.$
\item If $d_1=d_2-1$, then
$\nu(C)= 2(d_3-d_1)+2 \geq 4.$
\item If $d_1=d_2$, then
$\nu(C)= 2(d_3-d_2)+3\geq 3.$
\end{enumerate}

\end{prop}

\begin{prop}
\label{propCB}
Suppose that $C$ has type $2B$ and  exponents $(d_1,d_2,d_3,d_4)$.
Then the following hold. 
\begin{enumerate}
\item If $d_1 <d_2-2$, then
$\nu(C)=d_3+d_4-2d_2+3\geq 3.$
\item If $d_1=d_2-2$, then
$\nu(C)= d_{3}+d_{4}-2d_{1}-1 \geq 3.$
\item If $d_1=d_2-1$, then
$\nu(C)= d_3+d_4-2d_1+1 \geq 3.$
\item If $d_1=d_2$, then
$\nu(C)= d_3+d_4-2d_2+2 \geq 2.$
\end{enumerate}

\end{prop}

The equality in Proposition \ref{propCB} (4) holds if and only if
$d_1=d_2=d_3=d_4$, see also \cite[Theorem 3.11 and Example 3.12]{3syz} in relation with such curves.

It is natural to ask the question whether in Corollary \ref{cor1}, in both cases, the two types $2A$ and $2B$ are possible for the curve $C'$. 
Our next result  gives a rather general construction to start with a plus-one generated line  arrangement $\A'$ and get  arrangements $\A$ of type $2B$, by adding a generic line. 

\begin{thm}
\label{thm0}
Let $\A'$ be a  line arrangement of type 1 with exponents $(m_1,m_2,m_3)$ with $m_1 \geq 3$. Let $L$ be a line such that the intersection $L \cap \A'$ consists of $\deg \A'=m_1+m_2$ points.
Then  the line arrangement $\A=\A' \cup L$ has type $2B$ and  the exponents of $\A$ are  $$(m_1+1,m_2+1,m_3+1, m_1+m_2-1).$$
\end{thm}

In the next result, we consider a line arrangement as a special type of conic-line arrangement, i.e. a conic-line arrangement where the number of conics is zero. The example of conic-line arrangement constructed in Theorem \ref{thm4} is very different from the conic-line arrangements $\CC$ constructed in the proof of Corollary \ref{cor11}, since here all the singularities situated on the conic are quasi-homogeneous, which is not the case for the singularity  $(\CC,(0:0:1))$ in general, see Remark \ref{rk11}.

\begin{thm}
\label{thm4}
Let $\CC'$ be a free conic-line arrangement with exponents $(m_1, m_2)$ with $m_1\geq 2$. Let $Q$ be a smooth conic such that the intersection $Q \cap \CC'$ consists of $2\deg \CC' = 2(m_1+m_2+1)$ points.
Then  the arrangement $\CC=\CC' \cup Q$ has type $2A$ with exponents  $(m_1+2,m_2+2,m_1+m_2+1)$.
\end{thm}

 \begin{rk}
\label{rkF1} 
It is interesting to note that the formula for $\tau(C)$ in Proposition \ref{propB} (1) can be obtained by setting $d_4=d_3+1$ in formula for $\tau(C)$ in Proposition \ref{propB} (2). Similarly, all the formulas for $\nu(C)$ in Proposition \ref{propCA} can be obtained by setting $d_4=d_3+1$ in the corresponding formulas for $\nu(C)$ in Proposition \ref{propCB}.
In particular, the sets of exponents
$(n_1+1,n_2+1,n_1+n_2-2)$ and $(n_1+1,n_2+1,n_1+n_2-2, n_1+n_2-1)$ yield the same Tjurina number $\tau(C)$ and the same freeness defect $\nu(C)$. For this reason the Tjurina number and the freeness defect cannot be used to show that the type $2B$ does not occur in Theorem \ref{thm4}. The same situation occurs below in Theorem \ref{thm2}. To discard the case $2B$ in both Theorems we need to use deep results from \cite{ST}.
 \end{rk}

We have seen that Theorem \ref{thm0} produces line arrangements of type $2B$.
To get families of line arrangements of type $2A$,  one can use the following construction.
We start by constructing first a family of line arrangements $\A'$ of type $1$ by taking the union of two pencils of lines. The following result may be well-known to the specialists.
\begin{prop}
\label{prop2}
Let $\B_1$ (resp. $\B_2$) be the arrangement consisting  of $n_1$ lines passing through a point $a_1 \in \PP^2$ (resp. $n_2$ lines  passing through a point $a_2 \in \PP^2 \setminus \B_1$) and such that $1 \leq n_1 \leq n_2$ and $n_1+n_2>2$. Then the arrangement $\A'=\B_1 \cup \B_2$  has type $0$ and exponents $(1,n_1+n_2-2)$ for $n_1=1$, and type $1$ and exponents $(n_1,n_2,n_1+n_2-2)$ for $n_1 \geq 2$.
\end{prop}
In view of Theorem \ref{thm0}, this directly implies the following.
\begin{cor}
\label{cor2}
 If we add a generic line to the line arrangement $\A'$ from Proposition \ref{prop2} with $n_1 \geq 3$, then we get a line arrangement of type $2B$ with exponents $$(n_1+1, n_2+1, n_1+n_2-1, n_1+n_2-1).$$
\end{cor}
The next result tells us what happens when we add a line $L$ to the arrangement $\A'$ in Proposition \ref{prop2}, passing through exactly one double point of $\A'$.

\begin{thm}
\label{thm2}
Let $\A'$ be the line arrangement from Proposition \ref{prop2} and assume that  $3 \leq n_1 \leq n_2$. Let $L$ be a line such that the intersection $L \cap \A'$ consists of $n_1+n_2-1$ points.
Then  the arrangement $\A=\A' \cup L$ has  has type $2A$ with exponents $(n_1+1,n_2+1,n_1+n_2-2)$.
\end{thm}

Finally, we present an enumerative description of  line arrangements of type $2$. Let us first fix the notation. Let $\mathcal{L} \subset \mathbb{P}^{2}$ be an arrangement of $d$ lines. We denote by $t_{k}:=t_{k}(\mathcal{L})$ the number of $k$-fold intersection points of $\mathcal{L}$, and we define the maximal multiplicity of $\mathcal{L}$ as follows:
$$m(\mathcal{L}) = {\rm max} \{k : \text{arrangement } \mathcal{L} \text{ admits } k\text{-fold points}\}.$$
\begin{thm}
\label{ll}
Let $\mathcal{L} \subset \mathbb{P}^{2}$ be an arrangement of $d$ lines of type $2$.  Then one has
$$m(\mathcal{L})\geq \bigg\lceil\dfrac{4d}{d+5}\bigg\rceil.$$
Moreover, the following hold.
\begin{enumerate}
    \item If $\mathcal{L}$ has type $2A$ and exponents $(d_1,d_2,d_3)$, then
    $$\sum_{r\geq 2}(r-1)t_{r} \geq d_{1}d_{2}+2.$$
    \item If $\mathcal{L}$ has type $2B$ and exponents $(d_1,d_2,d_3,d_4)$, then
    $$\sum_{r\geq 2}(r-1)t_{r} \geq d_{1}d_{2}+1.$$
\end{enumerate}
\end{thm}

Both lower bounds for $\sum_{r\geq 2}(r-1)t_{r} $ for  line arrangements of type $2$  given above are sharp and this  follows from Example \ref{eb} below.

The content of this paper is outlined below. Section 2 provides a recap of some necessary preliminaries. The subsequent sections then proceed to prove these results and provide a number of useful examples and related additional remarks.

\section{Some preliminaries} 

With the notation from Introduction, we have the following result.

\begin{thm}[{\cite[Theorem 1.2]{Dimca1}}]
\label{deff}
Let $C \, : f=0$ be a reduced plane curve of degree $d$ and $d_{1}= {\rm mdr}(C)$. Then the following hold.
\begin{enumerate}
    \item If $d_{1} \leq  (d-1)/2$, then $\nu(C) = (d-1)^{2}-d_{1}(d-1-d_{1})-\tau(C)$.
    \item If $d_{1} \geq (d-1)/2$, then
    $$\nu(C) = \bigg\lceil\frac{3}{4}(d-1)^{2}\bigg\rceil - \tau(C).$$
\end{enumerate}
\end{thm}

If we set $T=3(d-2)$, then the sequence $n(f)_k$ is symmetric with respect to the middle point $T/2$, that is one has
\begin{equation}
\label{E1}
n(f)_a=n(f)_b
\end{equation}
for any integers $a,b$ satisfying $a+b=T$, see \cite{Se}. It was shown in \cite[Corollary 4.3]{DPop} that the graded $S$-module  $N(f)$ satisfies a Lefschetz type property with respect to multiplication by generic linear forms. This implies, in particular, the inequalities:
\begin{equation}
\label{in} 
0 \leq n(f)_0 \leq n(f)_1 \leq \ldots \leq n(f)_{[T/2]} \geq n(f)_{[T/2]+1} \geq \ldots \geq n(f)_T \geq 0.
\end{equation}
Moreover, for a $3$-syzygy curve $C$ of degree $d$ with exponents $(d_1,d_2,d_3)$, we have the following formula for the initial degree of the graded module $N(f)$, see \cite[Theorem 3.9]{3syz}:
\begin{equation}
\label{ID}
\sigma (C)= \min \{k \ : \ N(f)_k \ne 0 \}=3(d-1)-(d_1+d_2+d_3).
\end{equation}

\subsection{Some key exact sequences }
\label{ESsec}

Let us come back to the setting of Theorem \ref{thm1}. To simplify the notation, in the sequel we set 
\begin{equation}
\label{deg}
e_1 = \deg f_1 = \deg C_1 \text{ and } e_2=\deg f_2 = \deg C_2.
\end{equation}
Assume in addition that $C_2$ is smooth and all the points in $C_1 \cap C_2$ are quasi-homogeneous singularities of $C=C_1 \cup C_2$.
Then there is an exact sequence  for any integer $k$ given by
\begin{equation}
\label{ES0}
0 \to E_{C_1}(k-e_2) \stackrel{f_2} \longrightarrow  E_C(k) \to i_{2*}\OO_{C_2}(-K_{C_2}-R+(k-1)D) \to 0,
\end{equation}
where $E_C= \widetilde{D_0(f)} $ (resp. $E_{C_1}= \widetilde{D_0(f_1)} $ is the sheafification
of the graded $S$-module $D_0(f)$ (resp. $D_0(f_1)$), which is a rank two vector bundle on $\PP^2$, $K_{C_2}$ is the canonical divisor on ${C_2}$, $R$ is the reduced scheme of $C_1 \cap C_2$ and the divisor $D'$ corresponds to the intersection of a line in $\PP^2$ with the curve $C_2$, see \cite{DIS,STY}. By taking the associated long cohomology exact sequence of the sequence \eqref{ES0}, we get
 \begin{equation}
\label{ES} 
0 \to D_0(f_1)_{k-e_2} \to D_0(f)_k \to H^0(C_2,\OO_{C_2}(D+(k-1)D')) \to 
\end{equation}
$$ \to N(f_1)_{k-e_2+e_1-1} \to N(f)_{k+e_1+e_2-1}   \to H^1(C_2,\OO_{C_2}(D+(k-1)D')) \to $$
$$\to H^2(\PP^2, E_{C_1}(k-e_2)) \to  H^2(\PP^2, E_{C}(k)) \to 0, $$
with $D = -K_{C_2}-R$, and hence  $\deg D=2-2g_2-r$,
 where $g_{C_2}$ is the genus of the smooth curve $C_2$ and $r$ is the number of points in the reduced scheme of $C_1 \cap C_2$. 
Moreover, one has
\begin{equation}
\label{ES100} 
\dim H^2(\PP^2, E_{C_1}(k-e_2))=\dim D_0(f_1)_{n_1+n_2-(k-e_2)-4} 
\end{equation}
and
$$
\dim H^2(\PP^2, E_{C}(k))=\dim D_0(f)_{(n_1+n_2+1)-k-4},
$$
see \cite{Dmax,Se}.

\begin{rk}
\label{rktheta}
If we identify $D_0(f)$ with the quotient $S$-module 
$$\overline D(f)=D(f)/SE,$$ then the morphism $\theta$ is just given by the multiplication by $f_1$.
Here $D(f)$ denotes the $S$-module of all derivations of $f$, that is
$$D(f)=\{(a,b,c) \in S^3 \ : \ af_x+bf_y+cf_z \in (f)\},$$
where $(f)$ is the principal ideal generated by $f$, and $SE$ denotes the submodule in $D(f)$ generated by the Euler derivation
$E=(x,y,z)$. The first non-zero morphism in \eqref{ES0} is multiplication by $f_2$ if we replace $D_0(f_1)$ and $D_0(f)$ by $\overline D(f_1)$ and $\overline D(f)$ respectively, see for instance \cite[Theorem 1.6]{STY}. This fact is used in the proof of Theorems \ref{thm4} and \ref{thm2} below.
\end{rk}

\subsection{Exact sequences  for line arrangements}
\label{ESLA}
Consider the triple of line arrangements $(\A',\A, \A'')$, where
$\A'=\A \setminus \{L\}$ with $L$ being a line in $\A$, and 
$\A''$ is the arrangement consisting of points on the line $L$ corresponding to the intersection set $\A' \cap L$.
By a change of coordinates we can take $L:x=0$, and we set
$$\A' :f' = 0, \  \A: f = 0 \text{ and } \A'':f'' = 0,$$
for some polynomial $f'' \in R=\C[y,z]$. Then there is an exact sequence of graded $S$-modules, sometimes called the Euler exact sequence of the triple $(\A',\A, \A'')$, given by
 \begin{equation}
\label{ES200} 
0 \to D(f')(-1) \to D(f) \to D(f''),
 \end{equation}
where $D(f')$, $D(f)$ and $D(f'')$ are as in Remark \ref{rktheta} and in addition $D(f'')$ is regarded as an $S$-module using the restriction of scalars given by the ring morphism $\psi: S \to R$, which replaces $x$ by $0$ in any polynomial $h \in S$. Here the first morphism
$D(f')(-1) \to D(f)$ is the multiplication by $x$, and the second morphism is given by 
$$\rho: D(f) \to D(f''),  \  \rho(a_1,a_2,a_3)= (\psi(a_2), \psi(a_3)),$$
see \cite[Proposition 4.45]{OT} and \cite[Lemma 2.6]{ST} that covers a more general case when $\A'$ is a conic-line arrangement, maybe under some quasi-homogeneity conditions. 
If we replace the modules $D(f')$, $D(f)$ and $D(f'')$ by the modules
$\overline D(f')$, $\overline D(f)$ and $\overline D(f'')$, the sequence 
obtained from \eqref{ES200} is no longer exact. However, if we pass to the associated coherent sheaves on $\PP^2$, we get the following short exact sequence
 \begin{equation}
\label{ES201} 
0 \to E_{\A'}(-1) \to E_{\A} \to i_*E_{\A''} \to 0,
\end{equation}
where we followed the notation from \eqref{ES0}, see \cite[Theorem 3.2]{Sch}. Using this exact sequence we can give the following useful
sufficient condition for the sequence \eqref{ES200} to be right-exact, i.e. for the morphism $\rho$ to be surjective.
\begin{thm}
\label{lower}
With the above notation, let 
$$
0 \rightarrow \oplus_{i=1}^{m-2} S[-c'_i]
\rightarrow \oplus_{i=1}^m S[-d_i]
\rightarrow D_0(\A') \rightarrow 0
$$
be a minimal free resolution of $D_0(\A')$. If $|\A''|>c_i'-1$ for all $i$, then
the Euler exact sequence 
$$
0 \rightarrow D(\A')[-1]
\stackrel{\cdot x}{\rightarrow} D(\A)
\stackrel{\rho}{\rightarrow} D(\A'') 
$$
is right exact.

\end{thm}
\proof 
Let $|\A''| = n$, and note that this implies $|\A'| \geq n$. 
Then the $S$-module $D(\A'')$ has two generators, one of degree $1$ and another of degree $n-1$. The first element is the Euler derivation, which is always in the image of $\rho$.
It follows that 
$$ D(\A)
\stackrel{\rho}{\rightarrow} D(\A'') $$
 is surjective if and only if the induced morphism
 $$\overline  D(\A)
\stackrel{\overline \rho}{\rightarrow} \overline D(\A'') $$
is surjective.
Moreover, it suffices to show that $ \overline \rho$ is surjective in degree $n - 1$. Using the exact sequence \eqref{ES201}, it suffices to show that 
$H^1(\PP^2,E_{\A'}(n-2))=0$. 
By the minimal free resolution of $D_0(\A')$, we have 
\begin{eqnarray*}
0&=&H^1(\PP^2, \oplus\mathcal{O}(-d_i+n-2)) \rightarrow 
H^1(\PP^2,E_{\A'}(n-2)) \\
&\rightarrow &
H^2(\oplus\mathcal{O}(-c'_i+n-2)) \simeq H^0(\PP^2,\oplus \mathcal{O}(c'_i-n+2-3))=0.
\end{eqnarray*}
The last vanishing follows by our assumption, and the first one from the fact that $H^1(\PP^2, \mathcal{O}(k))=0$ for any integer $k$. This completes the proof.

\endproof
The exact sequence in Theorem \ref{lower} plays a key role in our proof of Theorem \ref{thm2}.
When we pass from line arrangements to conic-line arrangements, there are similar exact sequences according to \cite{ST}, and they are discussed in detail and used in the proof of Theorem \ref{thm4} below.

\section{The proofs of Theorems \ref{thm1} and \ref{thm10} and of Corollaries \ref{cor10} and \ref{cor11}} 

\subsection{Proof of Theorem \ref{thm1}}

Observe that there is a natural $S$-linear morphism
$$ \theta: D_0(f_1)_{k-e_2} \to D_0(f)_k$$
for any integer $k$. If $\delta \in D_0(f_1)_{k-e_2}$ is a syzygy, we can think about it as being a derivation killing $f_1$, and then we set
$$\theta (\delta)=f_2\delta -\frac{\delta(f_2)}{d}E,$$
where $d= \deg C=e_1+e_2$ and $E$ denotes the Euler derivation, see \cite[Theorem 5.1]{DIS}.
This morphism is injective. Indeed, if we suppose that $\theta (\delta)=0$ for some $\delta \in D_0(f_1)_{k-e_2}$, then in particular
$\theta (\delta)(f_1)=0$. Since $E(f_1)=e_1f_1 \ne 0$ and $\delta(f_1)=0$, we get $\delta(f_2)=0$. It follows that
$$0=\theta (\delta)=f_2\delta$$
and hence $\delta=0$ since $f_2 \ne 0$.

We denote by $d_1,d_2$ (resp. by $d_1', d_2'$) the first two exponents of $C_1$ (resp. of $C$). By Proposition \ref{propA}, we have
$$d_1' \leq d_1+e_2.$$
If $d_2' \leq d_2+e_2$, then by adding these two inequalities we get our claim. Suppose now that
$$d_2' > d_2+e_2.$$
 Let $r_1$ and $r_2$ be the first two generators in a minimal set of generators for $D_0(f_1)$ and let $r_1'$ be a generator of $D_0(f)$ of minimal degree $d_1'$. Then $\theta(r_1) \in D_0(f)_{d_1+e_2}$ and 
 $\theta(r_2 )\in D_0(f)_{d_2+e_2}$ are non-zero derivations. It follows there are non-zero homogeneous polynomials $p,q \in S$ such that $\theta(r_1)=pr_1'$ and $\theta(r_2)=qr_1'$.
 This implies $q\theta(r_1)-p\theta(r_2)=\theta(qr_1-pr_2)= 0$, and hence by the injectivity of $\theta$ (see the beginning of the proof of \cite[Theorem 5.1]{DIS}) we get
 $$qr_1-pr_2=0.$$
 This is a contradiction since the quotient $S$-module $D_0(f_1)/Sr_1$ is torsion free, as it follows, for instance, from the exact sequence in \cite[Theorem 5.1]{Split}. Indeed, this exact sequence may be rewritten as
 $$0 \to S(-d_1) \to D_0(f) \to B(C,r_1) \to 0,$$
 where the first morphism is given by $h \mapsto hr_1$. Hence
 $D_0(f_1)/Sr_1$ is isomorphic to an $S$-submodule $I$ of the ideal
 $B(C,r_1)$, in other words $I$ is an ideal in $S$, and hence has no torsion.
 
 \begin{rk}
\label{rk1}
Consider the Fermat curve $C_1: f_1=x^d+y^d+z^d=0$ and the line arrangement $C_2:f_2=x^d+y^d=0$. Then the curve $C_1$ has type
$d-1$ since the exponents of $C_1$ are $d_1=d_2=d_3=d-1$. On the other hand, the union $C=C_1 \cup C_2$ has type 0, since $C$ is a free curve, see \cite[Corollary 1.6]{DIPS}. It follows that there is no interesting  lower bound for the type of $C=C_1 \cup C_2$ without providing any additional information on curves $C_1$ and $C_2$.
\end{rk}

\subsection{Proof of Theorem \ref{thm10} and Corollaries \ref{cor10} and \ref{cor11}}

We use the exact sequence \eqref{ES} above, where
 $$g_2= \genus(C_2)=\frac{(e_2-1)(e_2-2)}{2}$$
 and $r=e_1e_2$.  Moreover, $D'$ consists of $e_2$ points, and hence
$$\deg (D+(k-1)D')=2-(e_2-1)(e_2-2)-e_1e_2+(k-1)e_2=e_2(k-e_1-e_2+2).$$
For $k < e_1+e_2-2$, we have $H^0(C_2,\OO_{C_2}(D+(k-1)D'))=0$, and this yields an isomorphism 
$$D_0(f_1)_{k-e_2} \to D_0(f)_k.$$
It follows that all the generators of $D_0(f)$ of degree $k < e_1+e_2-2$
come from generators of $D_0(f_1)$ of degree $k-e_2$.
In this way we get two generators of degrees
$$d_1' = d_1+e_2 \text{ and } d_2' = d_2+e_2.$$
A new generator, possibly coming from $H^0(C_2,\OO_{C_2}(D+(k-1)D'))=\C$, arises when $k=e_1+e_2-2$. Since
$$d_2'=d_2+e_2 \leq e_1+e_2 - 2$$
by our assumption, we get that $d_1'$ and $d_2'$ are the first two exponents of the curve $C$, and hence $C$ has type
$$d_1'+d_2'-(e_1+e_2-1)=(d_1+d_2-(e_1-1))+e_2=t_1+e_2$$
as claimed. 

The proof above implies that if the curve $C_1$ satisfies the condition
$d_1 \leq e_1-2$, then the resulting curve $C$ also satisfies this condition, since
$$d_2'=d_2+e_2 \leq e_1+e_2-2.$$
This completes the proof of Theorem \ref{thm10}.

\begin{ex}
\label{ex10}
The condition $d_2 \leq \deg(C_1)-2$ is necessary in Theorem \ref{thm10} as the following example shows. Let $C_1:f_1=x^2+y^2+z^2=0$ and $C_2:f_2=x=0$. Then $d_1=d_2=1$ and hence 
$t(C_1)=1$. A direct computation shows that $C=C_1 \cup C_2$ is a 3-syzygy curve with exponents $(1,2,2)$ and hence of type
$$1+2-2=1<t(C_1)+\deg(C_2)=2.$$
\end{ex}
\noindent
Now we prove Corollary \ref{cor10}. 
\begin{proof}
The first claim follows from
\eqref{bound}, since
$$t(\A)=d_1+d_2-(d-1) \leq 2d-4-(d-1)=d-3.$$
Moreover, \eqref{bound} tells us that we can apply Theorem \ref{thm10} to any line arrangement $C_1$.
Starting with the integers $t$ and $d$ as in the second claim of Corollary \ref{cor10}, we choose $C_1$ to be a free line arrangement of $e_1=d-t \geq 3$ lines, for instance $x^{e_1}+y^{e_1}=0$. Then we add $t$ generic lines to this arrangement and we get by Theorem \ref{thm10}
a line arrangement $\A$ of degree $d$ and type $t$.
This completes the proof of Corollary \ref{cor10}.
\end{proof}
Now we present our proof of Corollary \ref{cor11}.
\begin{proof}
We first notice that our additional assumption on $\CC$ implies
$$t(\CC)=d_1+d_2-\deg(\CC)+1 \leq 2d_2-\deg(\CC)+1\leq \deg(\CC)-3.$$
Then, for the second claim, it is enough as above to construct free conic-line arrangements $\CC_1$ for any $e_1=\deg(\CC_1)=d-t \geq 3$, containing at least one conic and such that $d_2 \leq \deg(\CC_1) - 2$.
When $e_1=3$, we take $\CC_1$ to be a smooth conic plus one of its tangents. It is known that then $\CC_1$ is free with exponents $(1,1)$ and we are done. 

For $e_1 \geq 4$, we consider the arrangement
$$\CC_1:f_1=x(xz+y^2)(x^{e_1-3}+y^{e_1-3})=0.$$
It is easy to check the following are two syzygies for $f_1$.
$$r_1=(x^2,xy,-e_1y^2-(e_1-1)xz)$$
and
$$r_2=(-(e_1-1)xy^{e_1-4},e_1x^{e_1-3}+y^{e_1-3},-2e_1x^{e_1-4}y+(e_1+1)y^{e_1-4}z).$$
It follows that 
$$d_1 \leq \deg(r_1)=2 \text{ and } d_2 \leq \deg(r_2)=e_1-3.$$
since  $r_1$ and $r_2$ are not  multiple of strictly lower degree syzygy. Using Theorem \ref{thmPO1}, it follows that $\CC_1$ is a free curve and the above inequalities are in fact equalities. In particular,
$$d_2=\max\{2,e_1-3\} \leq e_1-2$$
for $e_1 \geq 4$. 
To complete the proof of Corollary \ref{cor11} it remains to add $t$ generic lines to the arrangement $\CC_1$ and get in this way the required arrangement $\CC$ by Theorem \ref{thm10}.
\end{proof}
\begin{rk}
\label{rk11} 
The curve $\CC$ constructed above has a complicated singularity at the point $p=(0:0:1)$.
Indeed, looking at $(\CC,p)$ as the union of a $(e_1-2)$-ordinary point with the smooth germ at $p$ given by the conic $xz+y^2=0$ and using the known formula for the Milnor number of a union, see for instance \cite[Theorem 6.5.1]{CTC}, one can show that the Milnor number $\mu(\CC,p)$ is given by
$$\mu(\CC,p)=e_1^2-4e_1+6.$$

 On the other hand, using the formula for the total Tjurina number of a free curve, which can be obtained by setting $\nu(C)=0$ in Theorem \ref{deff} (1), we get
 $$\tau(\CC)=(e_1-1)^2-2(e_1-3)=e_1^2-4e_1+7.$$
 Now the curve $\CC$ has $e_1-3$ nodes $A_1$ as singularities except $p$, and hence
 $$\tau(\CC,p)=\tau(\CC)-(e_1-3)=e_1^2-5e_1+10.$$
 One can alternatively compute $\tau(\CC,p)$ by using only local computations, based on \cite[Proposition 3.2 and Theorem 5.1]{HH}.
 It follows that the invariant 
 $$\epsilon(\CC,p)=\mu(\CC,p)-\tau(\CC,p),$$
 which was introduced in \cite{POG} and measures how far $(\CC,p)$ is from being quasi-homogeneous, is given by
 $$\epsilon(\CC,p)=e_1-4.$$
 Hence for all $e_1 > 4$ the singularity $(\CC,p)$ is not quasi-homogeneous.
 
  \end{rk}
 
\section{The proof of Propositions \ref{propB}, \ref{propCA} and \ref{propCB}} 

\subsection{Proof of Proposition \ref{propB} }

To see that the number $m$ of syzygies can be only $3$ or $4$, and to get the corresponding values for the $\varepsilon_j$'s, use the equality \eqref{res2C}. In case (1), the value for $\tau(C)$ follows from
\cite[Proposition 2.1]{3syz} using that $d_1+d_2=d+1$, where $d= \deg C$. In case (2), we use \eqref{res2B} to get the following resolution of the Milnor algebra $M(f)$:
$$0 \to S(-d-d_3)\oplus S(-d-d_4) \to \oplus_{j=1}^{4}S(1-d-d_j) \to S^3(1-d) \to S.$$
This resolution implies that in this case one has
$$\tau(C)=\binom{k+2}{2}-3 \binom{k+3-d}{2}+\sum_{j=1}^{4}\binom{k+3-d-d_j}{2}-$$
$$-  \binom{k+2-d-d_3}{2}-  \binom{k+2-d-d_4}{2},$$
for any large $k$. The last claim follows from the equality \eqref{res2C}.
\endproof

Here are some examples of curves $C$ having type $2$.

\begin{ex}
\label{extype2}

(i) Consider the Bolza curve $C : x^{5} -y^{2}z^{3}-xz^{4} = 0$. It is well-known that $C$ has a unique singularity of type $E_{8}$ and one can check that $C$ has type $2A$ with exponents $(2,4,4)$. 

(ii) Consider a conic-line arrangement $\mathcal{CL}_1\subset \mathbb{P}^{2}_{\mathbb{C}}$ given by
$$f_1 = (x^2 + y^2 - z^2 )(y-z)(x^2 - z^2 )=0.$$
This arrangement is known to be free and it has $3$ singularities of type $A_{1}$ and $3$ singularities of type $A_{3}$, so this gives us $\tau(\mathcal{CL}_1)=12$, see \cite[Example 2.9]{MP1}. If we add the conic $Q \, : y^2-xz=0$ to $\mathcal{CL}_1$, then the resulting arrangement $\mathcal{CL}$ has quasi-homogeneous singularities and has type $2B$ with exponents $(4,4,5,5)$. 

(iii) Consider a conic  arrangement $\mathcal{C}_1\subset \mathbb{P}^{2}_{\mathbb{C}}$ given by
$$f_1=  (-3x^2+xy+yz+xz)(-3y^2 + xy+yz+zx)(-3z^2 + xy+yz+zx)=0.$$
This arrangement is known to be free  with exponents $(2,3)$ and it has $3$ singularities of type $A_5$ and just one ordinary triple point $D_4$.
If we add the conic $Q:f_2=x^2+y^2+z^2=0$, we get new conic  arrangement $\mathcal{C}$ which has type $2A$ and exponents
$(4,5,6)$. Note that the conic $Q$ is generic with respect to the arrangement $\mathcal{C}_1$, namely the intersection $\mathcal{C}_1\cap Q$ consists of $12$ points, which give rise to $12$ nodes for the arrangement $\mathcal{C}$.
\end{ex}

\begin{rk}
\label{rk2AB}
One very natural way to obtain line arrangements of type $2$ is the following. We start with a free line arrangement $\A$ and let $\A_{i,j}$ be the line arrangement obtained from $\A$ by deleting two lines, namely
$L_i$ and $L_j$. Then \cite[Theorem 1.8]{Chu} gives sufficient condition such that $t(\A_{i,j})=2$. More precisely, with the notation from this result, we have the following possibilities.
The arrangements in \cite[Theorem 1.8]{Chu} (1) have type $2B$, and the same holds for the arrangements in case (2), when $c_1<c_2$. On the other hand, the arrangements in \cite[Theorem 1.8]{Chu} (2) have type $2A$ if $c_1=c_2$ and $t(\A_{i,j}) \ne1$. Moreover, in \cite[Example 4.5]{Chu}, the author produces a free line arrangement $\A$ of $11$ lines, such that the arrangement $\A_{3,4}$ has type $2A$ and the arrangements $\A_{2,11}$ and $\A_{3,5}$ have type $2B$. 
\end{rk}

\subsection{Proof of Proposition \ref{propCA} and \ref{propCB} }

We start with Proposition \ref{propCA}. Here we use Theorem \ref{deff} and Proposition \ref{propB}. Note that the condition $d_1 <d_2-2$ is exactly the condition $d_1<(d-1)/2$ which occurs in Theorem \ref{deff} (1). Hence, to prove the claim $(1)$ in Proposition \ref{propCA}, we use the formula for $\nu(C)$ given in Theorem \ref{deff} (1) and the formula for $\tau(C)$ given in Proposition
\ref{propB}, $(1)$. To prove the claims (2), (3) and (4), we use 
the formula for $\nu(C)$ given in
Theorem \ref{deff} (2) and the formula for $\tau(C)$ given in Proposition
\ref{propB} $(1)$. This completes the proof of Proposition \ref{propCA}.

The proof of Proposition \ref{propCB} is analogous.

\section{Proof of Theorem \ref{thm0}, Theorem \ref{thm4}, Proposition \ref{prop2} and Theorem \ref{thm2}} 

\subsection{Proof of Theorem \ref{thm0} }

We use the exact sequence \eqref{ES} above with $C_1=\A':f_1=0$ and $C_2=L:f_2=0$.
Hence $e_1=m_1+m_2=r$, $e_2=1$, $g_2= \genus(C_2)=0$. Moreover, $D'$ is just a point, and hence
$$\OO_{C_2}(D+(k-1)D')=\OO_L(k+1-m_1-m_2).$$
Let $(d_1,d_2, \ldots ,d_m)$ be the exponents of the curve $C=\A$.

First, consider the case $k <m_1+m_2-1=e_1-1$. Then $H^0(L,\OO_L(k+1-m_1-m_2))=0$ and we get an isomorphism
$$D_0(f_1)_{k-1} \to D_0(f)_k.$$
Note that both $k=m_1+1$ and $k=m_2+1$ satisfy the condition $k <m_1+m_2-1$, and hence we get that $d_1=m_1+1$ and $d_2=m_2+1$.
It follows that $d_1+d_2=\deg \A+1$, and hence $\A$ has type $2$.
To see what happens for the other exponents of $\A$, we look now at the case $k=e_1-1$. Since $\A'$ is a $3$-syzygy curve,
using \eqref{ID} we get
$$\sigma(\A')=\indeg N(f_1)=3(e_1-1)-(m_1+m_2+m_3)=2e_1-m_3-3.$$
Let $T_1=3(e_1-2)$ and recall that $\dim N(f_1)_j=\dim N(f_1)_{T_1-j}$
for any integer $j$ by \eqref{E1}. It follows that $N(f_1)_j=0$ for any
$j>T_1-\sigma(\A')=e_1+m_3-3$. Since $m_3 \leq e_1-2$ by \eqref{bound}, this yields an exact sequence
 \begin{equation}
\label{ES1} 
0 \to D_0(f_1)_{e_1-2} \to D_0(f)_{e_1-1} \to H^0(L,\OO_{L})  \to 0.
 \end{equation}
When $m_3=e_1-2$, then there are two new generators in $D_0(f)_{e_1-1}$,
one corresponding to the generator of $D_0(f_1)$ of degree $m_3$, and the second one corresponding to the generator
of $H^0(L,\OO_{L})=\C$. In this case $$(d_1,d_2,d_3,d_4)=(m_1+1, m_2+1, e_1-1, e_1-1)$$
as claimed. When $m_3 <e_1-2$, then the exact sequence \eqref{ES1} shows that there is exactly one new generator for $D_0(f)$ in degree $d_3=m_3+1 < e_1-1$, and then a new generator in degree $d_4 = e_1-1$ corresponding to the generator
of $H^0(L,\OO_{L})=\C$. \endproof

\subsection{Proof of Theorem \ref{thm4}}
 
 This rather long proof is divided into two steps to help the reader understand it more easily.
 
{\bf Step 1.}  We use again the exact sequence \eqref{ES} above with $C_1=\CC' : f_1=0$ and $C_2=Q : f_2=0$.
Hence we have $e_1=m_1+m_2+1$, $r=2e_1$, $e_2=2$, $g_2= \genus(C_2)=0$. Moreover, $D'$ consists of $2$ points and hence
$$\OO_{C_2}(D+(k-1)D')=\OO_L(2k-2e_1).$$
Let $(d_1,d_2, \ldots ,d_m)$ be the exponents of the curve $C=\CC$.

First, consider the case $k <e_1$. Then $H^0(L,\OO_L(2k-2e_1))=0$ and we get an isomorphism
$$D_0(f_1)_{k-1} \to D_0(f)_k.$$
Note that both $k=m_1+2$ and $k=m_2+2$ satisfy the condition $k < m_1+m_2+1$, and hence we get that $d_1=m_1+2$ and $d_2=m_2+2$.
It follows that $d_1+d_2=e_1+3=\deg \CC+1$, and hence $\CC$ has type $2$.
To see what happens for the other exponents of $\CC$, we look now at the case $k=e_1$. Since $\CC'$ is a free curve, one has $N(f_1)=0$ and this yields an exact sequence
 \begin{equation}
\label{ES10} 
0 \to D_0(f_1)_{e_1-2} \to D_0(f)_{e_1} \to H^0(L,\OO_{L}) \to 0.
 \end{equation}
This exact sequence  shows that there is exactly one new generator, call it $\al$, for $D_0(f)$ in degree $d_3=e_1$ corresponding to the generator of $H^0(L,\OO_{L})=\C$. 

{\bf Step 2.} 
If the conic-line arrangement $\CC$ has type $2A$, then the exponents of $\CC$ are $(d_1,d_2,d_3)=(m_1+2,m_2+2,m_1+m_2+1)$, and our result is proved. To show that this is the case, we have to show that there is no need for an additional generator for $D_0(f)$ beyond $\al$.

 A triple $(\CC',\CC,\CC'')$ of CL (short for conic-line) arrangements, where $\CC'= \CC \setminus \{Q\}$ with $Q \in \CC$ a smooth conic and $\CC''$  the 0-dimensional arrangement on $Q$ given by the reduced intersection $\CC' \cap Q$, is called quasi-homogeneous if all the singular points of $\CC'$ and of $\CC$ situated on $Q$ are quasi-homogeneous. 

 \begin{thm}
\label{thmES}  
If a triple $(\CC',\CC,\CC'')$ of CL is quasi-homogeneous, then the associated exact sequence \eqref{ES0} with $C_1=\CC'$, $C=\CC$ and $C_2=Q$ is exact.
 \end{thm}
 
 \proof
 This result is a stronger version of  \cite[Proposition 3.7]{ST}.
 The authors in \cite{ST} require in their result that  all the singular points of $\CC'$ and of $\CC$ are quasi-homogeneous, but this stronger assumption is not necessary. Indeed, as noted in  \cite{DIS} and in subsection \ref{ESsec} above, for a similar condition in \cite{STY}, this assumption is used in \cite{ST} only to compute the difference 
$$\tau(\CC)-\tau(\CC')=\deg J_f -\deg J_{f'},$$
where we set
$$\CC':f'=0  \text{ and } \CC:f=0,$$
see for instance the proof of \cite[Proposition 2.8 and Proposition 3.7]{ST}. The difference 
$\tau(\CC)-\tau(\CC')$ can be expressed in the following way
$$\sum_{p \in \CC''} (\tau(\CC,p)-\tau(\CC',p)),$$
since for a point $q \in \CC \setminus \CC''$ one has 
$\tau(\CC,q)-\tau(\CC',q)=0$, as the two singularities $(\CC,q)$ and $\tau(\CC',q)$ coincide for $q \notin Q$, $Q$ is smooth and we set by convention $\tau(\CC',p)=0$ when the germ $(\CC',p)$ is the empty germ. Lemma 2.7 in \cite{ST} is applied to the singularities $p \in \CC''$, since for quasi-homogeneous singularities the Milnor number coincide with the Tjurina number. That is why it is enough to require that the singularities of $\CC$ and $\CC'$ situated on $Q$ are quasi-homogeneous.
\endproof 

In \cite[Proposition 3.5]{ST}, the following exact sequence of $S$-modules is constructed
 \begin{equation}
\label{ES101} 
0 \to D(f')(-2) \to D(f) \to D(f'')_{\psi},
 \end{equation}
where $D(f)$, $D(f')$ are as in Remark \ref{rktheta}. To explain the equation $\CC'':f''=0$ and the $S$-module $D(f'')_{\psi}$, note that by a change of coordinates on $\PP^2$ we can assume that $Q:q=y^2-xz=0$.
Then we have an isomorphism
$$v: \PP^1 \to Q \text{ given by } v(s,t)=(s^2,st,t^2).$$
Using this isomorphism, we can regard the 0-dimensional arrangement $\CC''$ as a 0-dimensional arrangement in $\PP^1$ and hence it is defined by an equation $f''=0$, where $f'' \in R=\C[s,t]$ is some homogeneous polynomial. Hence the graded $R$-module $D(f'')$ is well-defined. Moreover, the map $v$ above induces a ring morphism
$$\psi=v^*: S \to R \text{ given by } v^*(g)=g(v(s,t)),$$
and the $S$-module $D(f'')_{\psi}$ is just the module obtained from the $R$-module $D(f'')$ by restriction of scalars using $\psi$.
Next, we consider the morphisms in the exact sequence \eqref{ES101}.
The first morphism $D(f')(-2) \to D(f) $ is given by the multiplication by $q$, the equation of $Q$, and hence it preserves the homogeneous components, that is 
$$q(D(f')(-2)_k)=q(D(f')_{k-2}) \subset D(f)_k.$$
The second morphism, call it $\rho:D(f) \to D(f'')_{\psi}$, is much more subtle.
If $(a_1,a_2,a_3) \in D(f)$, then it is shown that
$$\psi(a_1)=sQ_1, \  \psi(a_2)=\frac{tQ_1+sQ_2}{2} \text{ and } \psi(a_3)=tQ_2$$
for some polynomials $Q_1,Q_2 \in R$, see \cite{ST} for details. With this notation, one sets
$$\rho(a_1,a_2,a_3)=(Q_1,Q_2)$$
and hence 
$$\rho(D(f)_k) \subset D(f'')_{2k-1}.$$
After these preliminaries, we can continue our proof in Step 2.
In our situation, $\CC''$ consists of $2e_1$ points, hence $f''$ is a homogeneous polynomial of degree $2e_1$ without multiple factors.
The derivation $\al \in D_0(f)_{e_1}$ that was constructed at the end of Step 1 can be regarded as sitting in $\overline D(f)_{e_1}$ by
using the identification from Remark \ref{rktheta}.
Then any representative $\al_0 \in D(f)_{e_1}$ of $\al$ is clearly not in the image of the morphism
$D(f')(-2) \to D(f) $ and hence $\rho(\al_0) \ne 0$.
 If we assume $\al_0=(a_1,a_2,a_3)$, then 
  \begin{equation}
\label{K1} 
za_1-xa_3 \ne 0.
 \end{equation}
 Indeed, $za_1=xa_3$ implies that there is $h \in S$ such that
 $a_1=xh$ and $a_3=zh$. But this implies that
 $$\al_1=\al_0-hE=(0,a_2',0)$$
 would be a new representative of $\al$ and clearly $\rho(\al_1)=0$, a contradiction.
 
If $\CC$ has type $2B$, it means that there is an additional generator of degree $d_4 >d_3=e_1$. Using \eqref{boundQ}, we get
$d_4 \leq (e_1+2)-1$, and hence $d_4=e_1+1$. On the other hand, the exact sequence \eqref{ES10} which holds due to Theorem \ref{thmES}, where $e_1$ is replaced by $e_1+1$, shows that  $\overline D_0(f)_{e_1+1}$ modulo the image of $\overline D_0(f_1)_{e_1-1}$ is 2 dimensional. It follows that we need a new generator of degree $e_1+1$ if and only the vector space $V$ spanned by
$x \al_0, y \al_0$ and $z\al_0$ in $\overline D_0(f)_{e_1+1}$ modulo the image of $\overline D(f_1)_{e_1-1}$ is 1-dimensional.
Note that $\rho(E)=E''$, the Euler derivation corresponding to the ring $R$, and hence there is an induced morphism
$$\overline \rho: \overline D(f) \to \overline D(f'')_{\psi}.$$
The element
$$x\al_0=(xa_1,xa_2,xa_3)$$
is equivalent modulo $SE$ to
$$\be=x\al_0-a_1E=(0, xa_2-ya_1, xa_3-za_1).$$
Since $xa_3-za_1=s^2tQ_2-t^2sQ_1=st(sQ_2-tQ_1)$, it follows that
  \begin{equation}
\label{K2} 
\be_1=\overline \rho(x\al)= \overline \rho(\be) = (0,s(sQ_2-tQ_1)) \in \overline D(f'')_{\psi}.
 \end{equation}
 A similar computation shows that
  \begin{equation}
\label{K3} 
\be_2=\overline \rho(z\al)= (-t(sQ_2-tQ_1),0) \in \overline D(f'')_{\psi}.
 \end{equation}
 The elements $\be_1$ and $\be_2$ are linearly independent in
 $\overline D(f'')_{\psi}$. Indeed, no linear combination
 $c_1\be_1+c_2 \be_2$ with $c_1,c_2 \in \C$ can have the form
 $hE''=(sh,th)$ for some $h \in R$. Indeed
 $$c_1\be_1+c_2 \be_2=(-c_2t(sQ_2-tQ_1),c_1s(sQ_2-tQ_1))=(sh,th)$$
 would imply $(c_1s^2+c_2t^2)(sQ_2-tQ_1)=0$, which is a contradiction by \eqref{K1}. This shows that the classes of $x \al_0$ and of $z\al_0$ in $\overline D(f)$ are linearly independent as well. Therefore $\dim V \geq 2$ and this completes the proof of Theorem \ref{thm4}.

\begin{ex}
\label{ex4}
An illustration of Theorem \ref{thm4} occurs in  Example \ref{extype2} (iii), where $(m_1,m_2)=(2,3)$. Note that in Example
\ref{extype2} (ii) the condition that the intersection $Q \cap \A'$ consists of $2\deg \A'=2(m_1+m_2+1)=10$ points is not fulfilled, since
$\tau(\mathcal{CL})=24$ by Proposition \ref{propB},  and hence
$$\tau(\mathcal{CL}_1)+10 = 22\ne \tau(\mathcal{CL}).$$
This explains why the type $2B$ may occur here.
\end{ex}

\subsection{Proof of Proposition \ref{prop2} }

The claims  for $n_1=1$ and $n_1=2$ follow from \cite[Proposition 4.7 (4)]{DIM}. Consider now the case $n_1 \geq 3$.
Let $\B:g=0$ be the line arrangement obtained from $\A'$ by adding the line
$L'$ determined by the points $a_1$ and $a_2$. Then $\B$ is a supersolvable line arrangement, and $a_1$ is a modular point of multiplicity $m_1=n_1+1$.
Using \cite[Proposition 3.2]{JACO2019}, it follows that
$$mdr(g)= \min \{m_1-1, (n_1+n_2+1)-m_1\}=n_1.$$
Since any supersolvable arrangement is free, it follows that $\B$ is a free arrangement with exponents $(n_1,n_2)$. Next, we apply \cite[Theorem 1.3]{POG} to the curve $C=\A'$ and the line $L=L'$. In our case, $r=|\A' \cap L'|=2$ and $\epsilon :=\epsilon(\A,L')=0$ since all the singularities of line arrangements are quasi-homogeneous.
It follows that the curve $\A'$ is plus-one generated with exponents
$(n_1,n_2,n_1+n_2-2)$. Indeed, we have to use the Case $(3)$ in 
\cite[Theorem 1.3]{POG}, where $d=\deg(C)=\deg(\A')=n_1+n_2$.
In the notation of \cite[Theorem 1.3]{POG}, we have $d_1' = n_1$ and hence $d_1\geq d_1'-1=n_1-1$. It follows that $$d_1+1-\epsilon \geq n_1 \geq 3>r,$$
and hence by the last claim in \cite[Theorem 1.3]{POG} the curve $C$ is not free. Therefore we are in the case (3) of \cite[Theorem 1.3]{POG}, and this yields our claim.

\subsection{Proof of Theorem \ref{thm2} }

This proof is similar to our proof of Theorem \ref{thm4}, so we divide it into two steps.

{\bf Step 1.}
We use the exact sequence \eqref{ES} above with $C_1=\A':f_1=0$ and $C_2=L:f_2=0$.
Hence $e_1=n_1+n_2$, $e_2=1$, $g_2= \genus(C_2)=0$, $r=n_1+n_2-1$. Moreover, $D'$ is just a point and we have
$$\OO_{C_2}(D+(k-1)D')=\OO_L(k+2-n_1-n_2).$$
Let $(d_1,d_2, \ldots ,d_m)$ be the exponents of the curve $C=\A$.

Consider first the case $k <n_1+n_2-2=e_1-2$. Then $H^0(L,\OO_L(k+2-d_1-d_2))=0$ and we get an isomorphism
$$D_0(f_1)_{k-1} \to D_0(f)_k.$$
Note that for $n_1\geq 4$, both $k=n_1+1$ and $k=n_2+1$ satisfy the condition $k <n_1+n_2-2$, and hence we get that $d_1=n_1+1$ and $d_2=n_2+1$.
To see what happens for $n_1=3$, we look now at the case $k=e_1-2$. Since $\A'$ is a 3-syzygy curve,
using \eqref{ID} we get
$$\sigma(\A')=\indeg N(f_1)=3(e_1-1)-(n_1+n_2+n_1+n_2-2)=e_1-1.$$
Let $T_1=3(e_1-2)$ and recall that $\dim N(f_1)_j=\dim N(f_1)_{T_1-j}$
for any integer $j$ by \eqref{E1}. It follows that $N(f_1)_j=0$ for any
$j>T_1-\sigma(\A')=2e_1-5$. In particular, this yields an exact sequence
 \begin{equation}
\label{ES11} 
0 \to D_0(f_1)_{e_1-3} \to D_0(f)_{e_1-2} \to H^0(L,\OO_{L})  \to 0.
 \end{equation}
Consider first the case $n_1=3<n_2$, when the sequence becomes
$$0 \to D_0(f_1)_{n_2} \to D_0(f)_{n_2+1} \to H^0(L,\OO_{L})  \to 0.$$
This shows that there are two new generators in $D_0(f)_{n_2+1}$,
one corresponding to the generator of $D_0(f_1)$ of degree $n_2$ and the other one, call it $\al$,  corresponding to the generator
of $H^0(L,\OO_{L})=\C$. This implies that $d_2=d_3=n_2+1$ in this case. When $n_1=3=n_2$, then there are three new generators in $D_0(f)_{n_2+1}$,
two corresponding to the two generators of $D_0(f_1)$ of degree $n_1=n_2$, and the third one, call it again $\al$,  corresponding to the generator
of $H^0(L,\OO_{L})=\C$. In this case $(d_1,d_2,d_3)=(4,4,4)$
as claimed.
It follows that in all cases $\A$ has type $2$.
Now consider the case $n_1 \geq 4$ and study the exponents $d_3$ and $d_4$.
The exact sequence \eqref{ES11} shows that there is exactly one new generator for $D_0(f)$ in degree $d_3=e_1-2=n_1+n_2-2$, which we denote by $\al$, corresponding to the generator
of $H^0(L,\OO_{L})=\C$. 

{\bf Step 2.}
If $\A$ is a $3$-syzygy curve, then Theorem \ref{thm2}  is proved. We show now that this is indeed the case.

The arrangement $\A''$ consists of $n=e_1-1=n_1+n_2$ points, using the formulas \eqref{res2B} and \eqref{res2C}, we see that 
$m=3$, $\varepsilon_1=1$ and $c_1=2(n_1+n_2-1).$
Hence the integer $c'_1$ from Theorem \ref{lower} is given by
$$c_1'=c_1-(e_1-1)=n_1+n_2-1.$$
It follows that we can apply Theorem \ref{lower} and get a surjective morphism $\rho$.
Now our claim that $D_0(f)$ can be generated by 3 elements or, equivalently, that $D(f)$ is generated by 4 elements, follows from the next basic fact in commutative algebra.
If we have an exact sequence of $S$-modules
$$0 \to M' \to M \to M'' \to 0,$$
$A$ a generating set for $M'$ and a set $B \subset M$ which maps onto a generating set for $M''$, then the image of $A$ under $M' \to M$ union $B$ gives a generating set for $M$.
We apply this to $M'=D(f')$, $M=D(f)$ and $M''=D(f'')$, and $A$ consisting of $E'=E$ the Euler derivation of $\A'$ and of $\A$, plus the  two generators of $D(f')_0$ and $B$ consisting of
$E$ and an element $\al \in D(f)$ such that $\rho(\al)$ is the generator of $D(f'')$ of degree $n-1$. The element $xE$ is discarded from the union since we have $E$ in this union.
\endproof

\section{The proof of Theorem \ref{ll}}

Let us recall that by \cite[Theorem 2.1]{DS14} if $C$ is a reduced curve of degree $d$ in $\mathbb{P}^{2}$ having only quasi-homogeneous singularities, then $$d_{1} \geq \alpha_{C}\cdot d - 2,$$
where $\alpha_{C}$ denotes the Arnold exponent of $C$. Let us recall that the Arnold exponent of a given curve $C$ is the minimum over all log canonical thresholds of singular points of $C$. Coming back to our setting, since for line arrangements all singular points are quasi-homogeneous, we can use the aforementioned result, but first we need to find $\alpha_{\mathcal{L}}$. Using \cite[Theorem 1.3]{Cheltsov}, we see that
$$\alpha_\mathcal{L} = \frac{2}{m(\mathcal{L})},$$
and we get
$$d_{1} \geq \frac{2d}{m(\mathcal{L})}-2.$$
Since $d_{1} \leq (d+1)/2$, we finally obtain
$$m(\mathcal{L}) \geq \bigg\lceil\dfrac{4d}{d+5}\bigg\rceil.$$

To prove the claim in (1), we use Proposition \ref{propB} (1).
Since $d^{2}=d_{1}^{2}+d_{2}^{2}+2d_{1}d_{2}-2d_{1}-2d_{2}+1$, we get
$$\tau(\mathcal{L}) = d^{2}-d_{1}d_{2}-2d_{3}-1.$$
Recall that 
$$\tau(\mathcal{L}) = \sum_{r\geq 2}(r-1)^{2}t_{r} = \sum_{r\geq 2}(r^{2}-2r+1)t_{r}.$$
Since $d^{2}-d = \sum_{r\geq 2}(r^{2}-r)t_{r}$, we have
$$\tau(\mathcal{L}) = d^{2}-d + \sum_{r\geq 2}(-r+1)t_{r}.$$
This gives us
$$d^{2}-d + \sum_{r\geq 2}(-r+1)t_{r} = d^{2}-d_{1}d_{2}-2d_{3}-1,$$
and finally
$$\sum_{r\geq 2}(r-1)t_{r} = d_{1}d_{2}-d+2d_{3}+1 = d_{1}d_{2} - d_{1} - d_{2} + 2d_{3}+2.$$
The claim follows now from the fact that $d_{1}\leq d_{2}\leq d_{3}$.

To prove the claim in (2), we use Proposition \ref{propB} (2). As above we get
$$\tau(\mathcal{L}) = d^{2}-d_{1}d_{2}-d_{3}-d_4.$$
This yields
$$d^{2}-d + \sum_{r\geq 2}(-r+1)t_{r} = d^{2}-d_{1}d_{2}-d_{3}-d_4,$$
and hence
$$\sum_{r\geq 2}(r-1)t_{r}=d_1d_2+d_3+d_4-d_1-d_2+1.$$
The claim (2) follows, since $d_1 \leq d_2 \leq d_3 \leq d_4$.

\begin{ex}
\label{eb}
We start by recalling \cite[Example 1.1]{3syz}. Let us consider the line arrangement $\mathcal{L} \subset \mathbb{P}^{2}$ given by
$$Q(x,y,z) = xyz(x+y-2z)(x-3y+z)(-5x+y+z)(x+y+z).$$
This arrangement has $d=7$ with $t_{3}=3$ and $t_{2} = 12$. One can check that $\mathcal{L}$ has type $2A$ with exponents $(4,4,4)$ and $\nu(\mathcal{L})=3$.

Moreover, a generic arrangement of $5$ lines described in Example \ref{extype2} has $10$ nodes and exponents $(3,3,3,3)$. For this arrangement of type $2B$ we have equality in Theorem \ref{ll} (2). 
\end{ex}

\section*{Acknowledgement}
We would like to thank the referees for their careful reading of the manuscript and for their very useful suggestions to improve the presentation.

Takuro Abe is partially supported by JSPS KAKENHI Grant Numbers JP23K17298 and JP23K20788.

Alexandru Dimca is partially  supported by the project ``Singularities and Applications'' - CF 132/31.07.2023 funded by the European Union - NextGenerationEU - through Romania's National Recovery and Resilience Plan.

Piotr Pokora is supported by the National Science Centre (Poland) Sonata Bis Grant \textbf{2023/50/E/ST1/00025}. For the purpose of Open Access, the author has applied a CC-BY public copyright license to any Author Accepted Manuscript (AAM) version arising from this submission.

\section*{Conflict of Interests}
We declare that there is no conflict of interest regarding the publication of this paper.

\section*{Data Availability Statement}
We do not analyze or generate any datasets, because this work proceeds within a theoretical and mathematical approach.

\end{document}